\documentclass[11pt]{article}

\usepackage{amsmath}
\usepackage{amsfonts}
\usepackage{amssymb}
\usepackage{graphicx}
\usepackage{epsf}
\usepackage{epsfig}
\usepackage[usenames]{color}
\usepackage[applemac] {inputenc}
\usepackage[T1]{fontenc}

\usepackage{fullpage}
\def\RR{\rm \hbox{I\kern-.2em\hbox{R}}}
\def\NN{\rm \hbox{I\kern-.2em\hbox{N}}}
\def\ZZ{\rm {{\rm Z}\kern-.28em{\rm Z}}}
\def\CC{\rm \hbox{C\kern -.5em {\raise .32ex \hbox{$\scriptscriptstyle
|$}}\kern
-.22em{\raise .6ex \hbox{$\scriptscriptstyle |$}}\kern .4em}}
\def\vp{\varphi}
\def\<{\langle}
\def\>{\rangle}
\def\t{\tilde}

\def\e{\varepsilon}

\def\nl{\newline}

\def\cN{{\cal N}}

\def\cF{{\cal F}}

\def\cE{{\cal E}}

\def\cD{{\cal D}}

\def\cO{{\cal O}}
\def\R{\mathbb{R}}
\def\N{\mathbb{N}}

\def\Chi{\raise .3ex
\hbox{\large $\chi$}} \def\vp{\varphi}
\def\lsima{\hbox{\kern -.6em\raisebox{-1ex}{$~\stackrel{\textstyle<}{\sim}~$}}\kern -.4em}
\def\lsim{\hbox{\kern -.2em\raisebox{-1ex}{$~\stackrel{\textstyle<}{\sim}~$}}\kern -.2em}
\def\gsim{\hbox{\kern -.2em\raisebox{-1ex}{$~\stackrel{\textstyle>}{\sim}~$}}\kern -.2em}
\def\[{\Bigl [}
\def\]{\Bigr ]}
\def\({\Bigl (}
\def\){\Bigr )}
\def\[{\Bigl [}
\def\]{\Bigr ]}
\def\({\Bigl (}
\def\){\Bigr )}

\newcommand{\be}{\begin{equation}}
\newcommand{\ee}{\end{equation}}
\newcommand{\bea}{$$ \begin{array}{lll}}
\newcommand{\eea}{\end{array} $$}
\newcommand{\bi}{\begin{itemize}}
\newcommand{\ei}{\end{itemize}}
\newcommand{\iref}[1]{(\ref{#1})}
\newtheorem{theorem}{Theorem}[section]
\newtheorem{remark}{Remark}[section]
\newtheorem{lemma}{Lemma}[section]

\newtheorem{corollary}{Corollary}[section]

\def\E{\mathbb E}

\def\N{\mathbb N}

\newcommand{\M}[1]{{\rm\textbf{M}}$_{#1}$}

\begin{document}

\title
{Sparse polynomial approximation of
parametric elliptic PDEs \\
Part II: lognormal coefficients
\thanks{ 
   This research was supported by the ONR Contracts
  N00014-11-1-0712,  N00014-12-1-0561, N00014-15-1-2181; the  NSF Grants  DMS 1222715;
  the Institut Universitaire de France; and  the ERC AdG BREAD. 
}
}
\author{ 
Markus Bachmayr, Albert Cohen, Ronald DeVore and Giovanni Migliorati 
    }

\maketitle
\date{}
\begin{abstract}
We consider the linear elliptic equation $-{\rm div}(a\nabla u)=f$ on some bounded domain $D$, where 
$a$ has the form $a=\exp(b)$ with $b$ a random function defined 
as $b(y)=\sum_{j\geq 1} y_j\psi_j$
where $y=(y_j)\in \R^\N$ are i.i.d.\ standard scalar Gaussian variables and $(\psi_j)_{j\geq 1}$ is a given sequence
of functions in $L^\infty(D)$.
We study the summability properties of
Hermite-type expansions of the solution map $y\mapsto u(y) \in V:= H^1_0(D)$, that is, expansions
of the form $u(y)=\sum_{\nu\in\cF} u_\nu H_\nu(y)$, where $H_\nu(y)=\prod_{j\geq 1} H_{\nu_j}(y_j)$ are the tensorized Hermite polynomials
indexed by the set $\cF$ of finitely supported sequences of nonnegative integers.
Previous results \cite{HS} 
have demonstrated that, for any $0<p\leq 1$, the 
$\ell^p$ summability of the sequence $(j\|\psi_j\|_{L^\infty})_{j\geq 1}$
implies $\ell^p$ summability of the sequence $(\|u_\nu\|_V)_{\nu\in\cF}$.
Such results ensure convergence rates $n^{-s}$ with $s=\frac 1 p-\frac 1 2$ of polynomial approximations
obtained by best $n$-term truncation of Hermite series, where the error is measured in
the mean-square sense, that is, in $L^2(\R^\N,V,\gamma)$, where $\gamma$ is the infinite-dimensional Gaussian measure.
In this paper we considerably improve these results by providing sufficient conditions for the
$\ell^p$ summability of $(\|u_\nu\|_V)_{\nu\in\cF}$ expressed in terms of the
pointwise summability properties of the sequence $(|\psi_j|)_{j\geq 1}$. This leads to a refined
analysis  which takes into account the amount of overlap between the supports
of the $\psi_j$. For instance, in the case of disjoint supports, our results
imply that, for all $0<p<2$ the $\ell^p$ summability of $(\|u_\nu\|_V)_{\nu\in \cF}$
follows from  the weaker assumption that $(\|\psi_j\|_{L^\infty})_{j\geq 1}$ is $\ell^q$ summable for
$q:=\frac {2p}{2-p}>p$. In the case of arbitrary supports, our results imply that 
 the $\ell^p$ summability of $(\|u_\nu\|_V)_{\nu\in \cF}$
follows from the $\ell^p$ summability of $(j^\beta\|\psi_j\|_{L^\infty})_{j\geq 1}$ for some $\beta>\frac 1 2$, which still represents
an improvement over the condition in \cite{HS}. We also explore intermediate cases
of functions with local yet overlapping supports, such as wavelet bases. One interesting 
observation following from our analysis is that for certain relevant examples, the use of the Karhunen-Lo\`eve
basis for the representation of $b$ might be suboptimal compared to other representations,
in terms of the resulting summability properties of $(\|u_\nu\|_V)_{\nu\in \cF}$.
While we focus on the diffusion equation, our analysis applies to other type of linear PDEs with 
similar lognormal dependence in the coefficients.
\end{abstract}

\noindent
{\bf Keywords:} stochastic PDEs, lognormal coefficients, $n$-term approximation, Hermite polynomials.
\vskip .1in

\section{Introduction}

\subsection{Approximation of high-dimensional parametric PDEs}

Parametric partial differential equations equations have
the general form
\be
\cD(u,y)=0,
\label{ppde}
\ee
where $u\mapsto \cD(u,y)$ is a partial differential operator
that depends on a parameter vector
$y=(y_j)_{j=1,\dots,J}$ ranging in a certain domain $U\in\R^J$, where $J\geq 1$ is fixed. 
Assuming well-posedness of the problem in some Banach space $V$,
the solution map 
\be
y \mapsto u(y),
\label{map}
\ee
is defined from the parameter domain $U$ to the solution space $V$.

Equations of this type arise both in stochastic and deterministic
modeling, depending on the nature of the parameters $y_j$
which may either be random or deterministic variables. In both settings,
one main computational challenge is to approximate the entire solution map $y\mapsto u(y)$
up to a prescribed accuracy, with reasonable computational cost. This task
has been intensively studied since the 1990s, see in particular 
\cite{GS1,GS2,KL,X} for general treatments.
It becomes very challenging when the number of parameters $J$ is large
due to the {\it curse of dimensionality}. Ideally, one would like to design numerical
methods that are immune to the growth of $J$, which in principle amounts to 
treat the case of {\it countably} many variables, that is,
\be
y=(y_j)_{j\geq 1}\in U \subset \R^\N.
\ee
This problem has been the object of much attention in recent years \cite{BNTT1,BNTT2,CD,CDS,CCS,G,GS,HS}. 

Sparse polynomial methods are based on approximations to $u$ of the form
\be
u_\Lambda(y):=\sum_{\nu\in \Lambda} u_\nu y^\nu,
\label{sparselambda}
\ee
where $\Lambda\subset \cF$ is a finite set of (multi-)indices $\nu=(\nu_j)_{j\geq 1} \in\cF$
and  $y^\nu=\prod_{j\geq 1} y_j^{\nu_j}$.
In the case of an infinite number of parameters, the index set $\cF$ 
denotes the (countable) set of 
all sequences of nonnegative integers which are 
{\em finitely supported} (i.e. those sequences for which
only finitely many terms are nonzero). Note that the polynomial
coefficients $u_\nu$ are functions in $V$,
and therefore the construction of $u_\Lambda$ requires in principle
the computation of $\#(\Lambda)$ such functions. 

One particularly
relevant example is the model elliptic PDE 
\be\label{pde}
-{\rm div}(a\nabla u)=f,
\ee
set on a bounded Lipschitz domain $D\subset \R^d$ with homogeneous Dirichlet boundary conditions (in typical applications 
$d=1,2,3$), where $a=a(y)$ is a diffusion coefficient that depends on $y$, and where $f\in H^{-1}(D)$ is fixed. The so-called {\it affine} case refers to a 
diffusion coefficient of the form
\be
a(y)=\bar a+ \sum_{j\geq 1} y_j\psi_j,
\ee
where $\bar a$ and the $(\psi_j)_{j\geq 1}$ are given functions from $L^\infty(D)$. In this case, the $y_j$ typically 
range on finite intervals which upon  renormalization of $\bar a$ and $\psi_j$ can be assumed to be $[-1,1]$. Therefore
the parameter domain is the infinite-dimensional box
\be
U:=[-1,1]^\N,
\ee
and well-posedness of the problem in $V=H^1_0(D)$ is ensured
for all $y\in U$ by the so-called {\em uniform ellipticity assumption}
\be
\sum_{j\geq 1}|\psi_j(x)|\leq \bar a(x)-r, \quad x\in D,
\ee
for some fixed $r>0$. For this model problem, convergence results have been obtained for polynomial 
approximations $u_{\Lambda_n}$, where $\#(\Lambda_n)=n$,  constructed by best $n$-term truncation of infinite polynomial 
expansions, either of Taylor or Legendre type, that is, by retaining the $n$ coefficients of largest norms
in such expansions.

A striking result, first established in \cite{CDS} under the uniform ellipticity assumption, states that whenever
  the sequence $(\|\psi_j\|_{L^\infty})_{j\geq 1}$ is $\ell^p$ summable for some $0<p<1$,
then such $n$ term polynomial approximations converge with rate $n^{-s}$ in $L^\infty(U,V)$
where $s:=\frac 1 p-1$. For the Legendre approximations, an  improved rate $n^{-s}$, with 
$s:=\frac 1 p-\frac 1 2$, is achieved in $L^2(U,V,\mu)$ where $\mu$ is the multivariate
uniform measure, that is, in the mean-square sense if the $y_j$ are i.i.d.\ uniformly distributed
variables. These results have been extended to a large range of linear or nonlinear
parametric PDEs \cite{CD,CCS} where $y$ is again ranging in the infinite dimensional box $U$.
They heavily rely on the holomorphy of the solution map $y\mapsto u(y)$ in each
variable $y_j$.

\subsection{Elliptic PDEs with lognormal coefficients}

In the present paper, we focus our attention 
on the so-called {\it lognormal} case for the elliptic PDE \iref{pde}. In this case,
the diffusion coeffient $a$ is of the form
\be
a=\exp(b),
\ee
where $b$ is a random function of
the form
\be
b=b(y)=\sum_{j\geq 1}y_j\psi_j,
\label{paramb}
\ee 
defined from a given sequence
$(\psi_j)_{j\geq 1}$ of functions in $L^\infty(D)$, with $y=(y_j)_{j\geq 1}$ a sequence of i.i.d.\ $\cN(0,1)$ variables. 
Thus, the parameter vector $y$ now ranges over the unbounded domain
\be
U:=\R^\N.
\ee
We work with the usual product measure space given by 
\be
(U, \mathcal{B}(U), \gamma)=(\R^\N, \mathcal{B}(\R^\N), \gamma),
\label{measure}
\ee 
where $\mathcal{B}(U)=\mathcal{B}(\R^\N)$ denotes
the $\Sigma$-algebra generated by the Borel cylinders
and $\gamma$ the tensorized Gaussian probability measure.

We discuss further several conditions on the family $(\psi_j)_{j\geq 1}$ which ensure that 
the series in \iref{paramb} converges almost surely in $L^\infty(D)$.
Under such conditions, $b(y)$ is a Gaussian random variable 
with values in $L^\infty(D)$. By the Lax-Milgram lemma, the associated 
weak solution $u(y)$ of \iref{pde} with $a=a(y)=\exp(b(y))$ is a random variable 
with values in $V=H^1_0(D)$.

One frequently used approach that leads to this framework is
by starting with $b=(b(x))_{x\in D}$ defined as a centered Gaussian process
over the domain $D$, with prescribed covariance function
\be
C_b(x,x'):=\E\bigl(b(x)b(x')\bigr),\quad x,x'\in D.
\ee 
One then obtains a representation of the form \iref{paramb} 
by considering the Karhunen-Lo\`eve decomposition
\be
b=b(\xi)=\sum_{j\geq 1} \xi_j \vp_j.
\ee
Here $(\vp_j)_{j\geq 1}$ are the $L^2$-orthonormal basis of eigenfunctions
of the integral operator with kernel $C_b$ 
and $\xi_j$ are independent centered Gaussian variables. One then 
sets
\be
y_j:=\lambda_j^{-1/2} \xi_j\quad {\rm and} \quad \psi_j:=\lambda_j^{1/2}\vp_j, \quad \lambda_j:=\E(\xi_j^2).
\ee
Here again, conditions on $(\psi_j)_{j\geq 1}$ are needed in order to ensure that
the series in \iref{paramb} converges almost surely in $L^\infty(D)$.

In this paper, we prefer to start with $b$ defined through \iref{paramb} for
a more general family $(\psi_j)_{j\geq 1}$. In particular we want to consider cases where
the family $(\psi_j)_{j\geq 1}$ in \iref{paramb} is nonorthogonal and therefore differs from
a normalized Karhunen-Lo\`eve basis of $b$. One typical example is the case where $b$ 
is a one-dimensional Brownian motion, which has a natural expansion in terms of the 
nonorthogonal Schauder basis. As shown later in this paper, this representation appears
to be more efficient than the Karhunen-Lo\`eve expansion in terms of the resulting best $n$-term
polynomial approximation rates.

In contrast to the affine model, the lognormal model is intrinsically stochastic: it reflects the situation where the diffusion
coefficient $a$ is allowed to become arbitrarily small or large (unlike in the 
affine case with the uniform ellipticity assumption), however with the probability 
being controlled by the Gaussian distribution. In this case, the relevant polynomial expansion
is the orthonormal Hermite series,
\be
u=\sum_{\nu\in\cF} u_\nu H_\nu, \quad   u_\nu :=\int_U u(y)H_\nu(y)\,d\gamma(y), \quad H_\nu(y)=\prod_{j\geq 1}H_{\nu_j}(y_j),
\label{multiherm}
\ee
where the univariate Hermite polynomials $H_k$ are normalized in $L^2(\R, dg)$, with $g$
denoting the standard Gaussian density. Note that 
 $(H_\nu)_{\nu\in\cF}$
is an orthonormal basis for the Hilbert space $L^2(U,\R,\gamma)$.
The coefficients $u_\nu$ are well defined in $V$ and the above series converges 
 in the Bochner space $L^2(U,V,\gamma)$ associated to $L^2(U,\R,\gamma)$
 whenever $u$ belongs to $L^2(U,V,\gamma)$.

Note that a function $u: y\mapsto u(y)$ in $L^2(U,V,\gamma)$ may also formally be viewed as a function 
$(x,y)\mapsto u(x,y):=u(y)(x)$ of both variables $x\in D$ and $y\in U$, and the inner product
between two functions $u,v\in L^2(U,V,\gamma)$ has the expression
\be
\<u,v\>=\int_U \int_D \nabla u(x,y)\cdot\nabla v(x,y)\,dx\, d\gamma(y),
\ee
where the operator $\nabla$ is always meant in the $x$ variable. In order to simplify notation, we shall
systematically avoid mentioning the $x$ variable, except when necessary, by simply writing the above integral as
\be
\<u,v\>=\int_U \int_D \nabla u(y)\cdot\nabla v(y) \,d\gamma(y).
\ee

The lognormal model and its approximation have been studied in various papers, e.g. \cite{BabNobTem07,Ch,DS,EMSU,G,GS,GKNSSS,HS,
KSSSU}.  Sufficient conditions for the finiteness of moments $\E(\|u(y)\|_V^k)$ at all orders $0\leq k<\infty$ have
been established, either by assuming smoothness properties
of the covariance kernel $C_b$ 
or by assuming a summability property for  the sequence  $(\|\psi_j\|_{L^\infty(D)})_{j\geq 1}$,
see for example \cite{Ch,DS,HS}. Since
\be
\E\bigl(\|u(y)\|_V^2\bigr)=\|u\|_{L^2(U,V,\gamma)}^2=\sum_{\nu\in\cF}\|u_\nu\|_V^2,
\ee
  the finiteness of  the second moment  is a necessary and sufficient condition for 
the $\ell^2$ summability of the  sequence $(\|u_\nu\|_V)_{\nu\in\cF}$ and the convergence of the Hermite
series in $L^2(U,V,\gamma)$. However, in order to  prove convergence
rates for the  best $n$-term truncation
\be
u_{\Lambda_n} :=\sum_{\nu\in\Lambda_n} u_\nu H_\nu,
\ee
of \iref{multiherm}, where $\Lambda_n$ denotes the set of indices corresponding to the $n$ largest $\|u_\nu\|_V$, we need to study the $\ell^p$ summability properties of
$(\|u_\nu\|_V)_{\nu\in\cF}$ for $p<2$. For example, if this sequence is proven to be $\ell^p$ summable, then a 
direct application of Stechkin's lemma \cite{De,CDS,CD}, 
implies
\be
\|u-u_{\Lambda_n}\|_{L^2(U,V,\gamma)}
=\(\sum_{\nu\notin\Lambda_n} \|u_\nu\|_V^2\)^{1/2} \leq C(n+1)^{-s}, 
\ee
where $s:=\frac 1 p-\frac 1 2$ and $C:=\bigl\| (\|u_\nu\|_V)_{\nu\in\cF} \bigr\|_{\ell^p(\cF)}$.  

The first and  currently the only available result concerning $\ell^p$ summability of Hermite coefficients, in
the infinite dimensional framework, has been established in \cite{HS} and reads as follows.

\begin{theorem}
\label{theohs} 
For any $0<p\leq 1$, 
if $(j\|\psi_j\|_{L^\infty})_{j\geq 1} \in \ell^p(\N)$ then $(\|u_\nu\|_V)_{\nu\in\cF}\in \ell^p(\cF)$.
\end{theorem}

The Hermite coefficients studied in \cite{HS} are actually different from those defined in \iref{multiherm},
since the authors are interested in a specific approximation process, namely the Galerkin
method. One specific difficulty in this approach is the fact that the bilinear form
\be
B(u,v)=\int_{U} \int_D a(y)\nabla u(y)\cdot \nabla v(y) \, d\gamma(y),
\ee
is not continuous and coercive on $L^2(U,V,\gamma)$. This particular issue has been dealt
with in papers such as \cite{Ch,GS,G,EMSU}. The approach from \cite{G}, which is also used in \cite{HS},
is based on a modified bilinear form where $\gamma$ is replaced by an auxiliary Gaussian measure.
The Galerkin error can then be estimated by the $\ell^2$ tails of modified Hermite coefficients corresponding to
a further auxiliary Gaussian measure. However, inspection of the proof of Theorem \ref{theohs} in \cite{HS} shows
that it also holds for the standard Hermite coefficients. Note that this issue is specific to the Galerkin
method, and does not occur in pseudo-spectral or least-squares methods.

In the present paper, we only consider the classical Hermite coefficients in \iref{multiherm}, 
since our main interest is to obtain general polynomial approximation
results for the solution map, that may further serve as a benchmark for any numerical method
based on such approximations.

The above Theorem \ref{theohs} is comparable to those obtained in \cite{CDS} for Taylor and Legendre
coefficients in the affine case, except for the appearance of the additional factor $j$ in front of $\|\psi_j\|_{L^\infty}$,
which makes the summability assumption more restrictive. On the one hand,
it is quite remarkable that such a result exists, since the  complex variable arguments used   in the affine case,  based 
solely on holomorphy of the solution map, are not sufficient to control the 
Hermite coefficients. The proof of the above result in \cite{HS} relies, instead, on estimates
for the mixed partial derivatives of $u$ up to some {\it finite} order $m$ related to $p$.
On the other hand, the type of assumptions used in this result
indicates that it is not in the sharpest possible form as we now explain.

First, the above result does not cover $\ell^p$ summability for $1<p<2$, which  one might expect to be obtainable
under weaker assumptions. Second, and more importantly, the summability conditions imposed on the $\|\psi_j\|_{L^\infty}$ in these results become
quite strong and artificial in the case when the supports of these functions do not overlap too much.

As a relevant example for this second point, consider the case where the $(\psi_j)_{j\geq 1}$ are a wavelet basis on the domain $D$. In this case
it is more natural to index such bases according to $(\psi_\lambda)_{\lambda\in\nabla}$, where $\lambda$ is 
a scale-space index following
the usual terminology, such as in \cite{Co}. With the notation $l=|\lambda|$ for the scale level, there
are $\cO(2^{dl})$ wavelets at this level and each of them has support of diameter $\cO(2^{-l})$. The supports
of wavelets at a given scale $l$ have finite overlap, in the sense that any $x\in D$ is contained in
the support of at most $M$ wavelets of level $l$ where $M$ is independent of $x$ and $l$. We assume 
that the $L^\infty$ norms of the wavelets only depend on the scale level, that is,
\be
\|\psi_\lambda\|_{L^\infty} = c_l, \quad |\lambda|=l.
\ee
It is well known that the geometric rate of decay of the wavelet contributions, as the scale level grows, reflects 
the amount of smoothness in the expansion (or the smoothness of the correlation function in the case of a random series).
It is thus natural to study the situation where $c_l$ is of the form
\be
\label{boundalpha}
c_l:=C2^{-\alpha l},
\ee
for some given $\alpha>0$.
Then, it can be checked, for example using the arguments from \cite{Ch},
which are recalled further in \S 2, 
that  for  arbitrary $\alpha>0$ and $C>0$,  the bound \iref{boundalpha} implies that
the solution map has bounded moments $\E\bigl(\|u(y)\|_V^k\bigr)$ for all $0\leq k<\infty$.  Indeed, this follows from the
Lax-Milgram lemma and the fact that
$\E\bigl(\exp(k\|b(y)\|_{L^\infty})\bigr)$ is finite. The particular case $k=2$ implies
that the Hermite coefficient sequence
$(\|u_\nu\|_V)_{\nu\in \cF}$ belongs to $\ell^2(\cF)$, provided only that $\alpha>0$. However, if we want to use the above mentioned result 
from \cite{HS} to prove
$\ell^p$ summability of this sequence for values of $p<2$, the appearance of the factor $j$, suddenly, requires  the  
very strong constraint $\alpha >2d$ which implies a strong smoothness condition on the diffusion coefficients.

The above wavelet example reveals a gap in the currently available analysis: $\ell^2$ summability can be obtained 
under mild assumptions on the smoothness of the diffusion coefficient, while proving $\ell^p$ summability for $p<2$
by the existing results immediatley imposes much higher smoothness (in the sense of the required decay of $\|\psi_j\|_{L^\infty}$ as $j\to \infty$).
The goal of the present paper is to propose a sharper analysis, which removes this gap. 
Let us mention that the same gap occurs, in a slightly less pronounced form, for the Legendre
coefficients in the affine case, an issue which we recently adressed in \cite{BCM}.

\subsection{Main results and outline of the paper}

The main result of our paper for the model elliptic equation \iref{pde} with lognormal diffusion 
coefficient is the following.

\begin{theorem}
\label{maintheo}
Let $p<2$ and let $q=q(p):=\frac {2p}{2-p}$. Assume that there exists a positive sequence $(\rho_j)_{j\geq 1}$ such that
\be
\sup_{x\in D}\sum_{j\geq 1} \rho_j|\psi_j(x)| <\infty.
\label{firstcond}
\ee
and
\be
(\rho_j^{-1})_{j\geq 1} \in \ell^q(\N).
\label{secondcond}
\ee
Then, the solution map $y\mapsto u(y)$ belongs to $L^k(U,V,\gamma)$ for all $0\leq k<\infty$. Moreover,
$(\|u_\nu\|_V)_{\nu\in \cF} \in \ell^p(\cF)$. In particular, best $n$-term Hermite approximations
converge in $L^2(U,V,\gamma)$ with rate $n^{-s}$ where $s=\frac 1 p-\frac 1 2=\frac 1 q$.
\end{theorem}

This theorem gives a significant improvement over Theorem \ref{theohs} in the
case where the functions $\psi_j$ do not overlap too much.   First of all, in the case of disjoint
supports, the first condition in Theorem \ref{maintheo} is met when $\rho_j^{-1}:=\|\psi_j\|_{L^\infty}$. Therefore,
it implies that, for all $0<p<2$ the $\ell^p$ summability of 
$(\|u_\nu\|_V)_{\nu\in \cF}$
follows from the assumption that $(\|\psi_j\|_{L^\infty})_{j\geq 1}$ is $\ell^q$ summable for
$q=q(p):=\frac {2p}{2-p}$. Note that $q(p)>p$ for any $p>0$ and that
\be
\lim_{p\to 2} q(p)= +\infty,
\ee
which shows that almost no decay of $(\|\psi_j\|_{L^\infty})_{j\geq 1}$
is required as $p$ gets closer to $2$.
Secondly, as is shown later in this paper, we can also use Theorem \ref{maintheo}  to treat the above mentioned
wavelet case, and obtain $\ell^p$ summability results with $p<2$ for {\it
any} smoothness index $\alpha>0$. Finally, even in the case of arbitrarily supported $\psi_j$, 
we   establish as a corollary of Theorem \ref{maintheo} that the $\ell^p$ summability of $(\|u_\nu\|_V)_{\nu\in \cF}$
follows from the $\ell^q$ summability of $(\|\psi_j\|_{L^\infty})_{j\geq 1}$ with $q:=\frac {2p}{2-p}$,
which again represents a significant improvement
over the condition in Theorem \ref{theohs}, since this $\ell^q$ summability follows by H\"older's  inequality from
the $\ell^p$ summability of $(j^\beta \|\psi_j\|_{L^\infty})_{j\geq 1}$ for some $\beta>\frac 1 2$.

Let us mention that a result analogous to Theorem \ref{maintheo} has been established 
in \cite{BCM} for Taylor and Legendre coefficients in the affine case, under the additional uniform ellipticity assumption.
A common feature with the present paper is that one key ingredient of  both proofs consists in establishing certain {\em weighted} $\ell^2$
estimates for the coefficients $\|u_\nu\|_V$ which further translate into $\ell^p$ estimates
by H\"older's inequality. However, establishing such weighted $\ell^2$ estimates for the Hermite
coefficients in the present lognormal cases uses completely different techniques than those
in \cite{BCM}. In particular, similar to \cite{HS}, we only rely on the study of mixed partial
derivatives $\partial^\mu u$ for a limited order $\|\mu\|_{\ell^\infty}\leq r$, while the estimates 
in \cite{BCM} exploit all orders $\mu\in\cF$. Let us also mention that,
while we focus on the diffusion equation, inspection of our proofs
reveals that our main results can be extended to other types of linear elliptic or parabolic PDEs with 
similar lognormal dependence of the coefficients.

The rest of our paper is organized as follows. We begin in \S 2 by revisiting conditions
which ensure that $y\mapsto u(y)$ is a measurable map with values in $V$ and
with finite moments $\E\bigl(\|u(y)\|_V^k\bigr)$ for all $0\leq k<\infty$. We introduce
a sufficient assumption in terms of the convergence in $L^\infty(D)$ of the series $\sum_{j\geq 1} \rho_j |\psi_j|$
for a positive sequence $(\rho_j)_{j\geq 1}$ which satisfies the summability property
\be
\sum_{j\geq 1} \exp(-\rho_j^2)<\infty.
\ee
This is a very weak condition on the sequence $(\rho_j)_{j\geq 1}$ and in particular always holds under the assumptions in Theorem \ref{maintheo}. 

As a first step in the proof of Theorem \ref{maintheo}, we relate in \S 3 the norms $\|u_\nu\|_V$
of Hermite coefficients with $\|\partial^\mu u\|_{L^2(U,V,\gamma)}$ 
for relevant values of $\mu$. This leads us to an identity between certain weighted
$\ell^2$ norms of both quantities. As mentioned above,
while all $\nu\in\cF$ are considered for the Hermite coefficients, only
limited order $\|\mu\|_{\ell^\infty}\leq r$ are considered for the partial derivatives.
As a second step, we obtain in \S 4 bounds on the previously introduced
weighted $\ell^2$ norms of the $\|\partial^\mu u\|_{L^2(U,V,\gamma)}$
under the assumption that $\sum_{j\geq 1} \rho_j |\psi_j(x)|$ is bounded
by a relevant constant. Finally, we combine these ingredients in \S 5 
to complete the proof of Theorem \ref{maintheo}.

In \S 6, we derive various consequences of Theorem \ref{maintheo}
corresponding to the different cases outlined above for the functions $(\psi_j)_{j\geq 1}$: (i) disjoint or
finitely overlapping supports, (ii) wavelets, (iii) arbitrary supports.
We discuss, in all three cases, the improvements over Theorem \ref{theohs}.

We conclude in \S 7 with an interesting observation 
which follows from our analysis, in the particular case where 
$b$ is a Brownian motion. For this case we compare the convergence
of the best $n$-term truncation of Hermite series when using
either the Karhunen-Lo\`eve representation  or the Schauder basis representation
for $b$. The first representation satisfies by construction an $L^2$ orthogonality between the $\psi_j$,
while the second one does not. Due to the limited amount of smoothness of $b$,
the assumptions of Theorem \ref{maintheo} do not hold for the Karhunen-Lo\`eve representation,
so that no algebraic convergence rate can be established for the resulting polynomial approximation
with our currently available techniques. In contrast, we can obtain an algebraic
convergence rate when using the Schauder representation, by exploiting
the local support properties of the basis functions. This hints that, in relevant practical cases, the 
Karhunen-Lo\`eve representation
might be suboptimal in terms of the resulting polynomial
approximation rates for the solution map  $y\mapsto u(y)$ in $L^2(U,V,\gamma)$.

Let us stress that our results only quantify the {\it approximability} of the solution map.
They should therefore be viewed as benchmark for concrete numerical methods.
The development of numerical methods that provably meet such benchmarks
will be the object of further investigation.

\section{Measurability and integrability of the solution map}

The first question we investigate is what  conditions on $(\psi_j)_{j\ge 1}$ guarantee that the solution map $y\mapsto u(y)$ is in $L^k(U,V,\gamma)$, i.e.
it is measurable with values in $V$ and satisfies
\be
\E\bigl(\|u(y)\|_V^k\bigr)=\int_U \|u(y)\|_V^k\, d\gamma(y) <\infty.
\label{squareint}
\ee
We recall that a function from a measurable space to a Banach space $B$ is measurable 
(also sometimes called $\mu$-measurable, or Bochner-measurable, or strongly measurable),
if and only if it is the almost everywhere pointwise limit of a sequence of simple functions.

Recall that $u(y)$ is the weak solution to the diffusion equation
\be
-{\rm div}(a(y)\nabla u(y))=f,
\ee
on the bounded domain $D$ with homogeneous Dirichlet boundary condition, for a fixed $f\in H^{-1}(D)$, were
\be
a(y):=\exp(b(y)).
\ee
The solution $u(y)$ is well defined as an element of $V=H^1_0(D)$ provided that $\sum_{j\geq 1} y_j\psi_j$ defines a function
$b(y)\in L^\infty(D)$. In such a case, by the Lax-Milgram lemma, we have
\be
\|u(y)\|_V \leq  C \|a(y)^{-1}\|_{L^\infty}\leq C\exp(\|b(y)\|_{L^\infty}),\quad C:=\|f\|_{H^{-1}}.
\label{laxmilgram}
\ee
This motivates the study of the finiteness of exponential moments  of $\|b(y)\|_{L^\infty}$ which we formulate for 
any given $0\leq k<\infty$ as follows.
\nl
\nl
\noindent
{\bf Property \M{k}:} The map $y\mapsto b(y)$ is measurable from $U$ to $L^\infty(D)$, and the exponential moment $\E\bigl(\exp(k\|b(y)\|_{L^\infty})\bigr)$
is finite.
\nl

According to \iref{laxmilgram}, we find that Property \M{k} implies the finiteness of the $k$-th moment $\E(\|u(y)\|_V^k)$.
Several types of conditions on $b$ have been introduced in 
the literature \cite{Ch,DS,HS} which are sufficient to guarantee that Property \M{k} holds.

One first approach requires that $b=(b(x))_{x\in D}$ is a Gaussian random field with some minimal smoothness,
in the sense of a smoothness assumptions on the covariance kernel
$C_b$. Indeed, if $C_b$ belongs to the H\"older space
$C^\beta(D\times D)$ for some $\beta>0$, that is, there exists $C>0$ such that
\be
|C_b(x,x')-C_b(z,z')| \leq C (|x-z|^\beta+|x'-z'|^\beta), \quad x,z,x',z' \in D,
\label{cond1}
\ee
then, following the argument of Proposition 2.1 in \cite{Ch}, we first find by using the Kolmogorov continuity theorem that there exists a version of 
$b$ with trajectories almost surely in $C^\alpha(D)$ for $\alpha<\beta/2$. By application of Fernique's theorem, see Proposition 2.3 of \cite{Ch}, we then find 
that 
\be
\E\bigl(\exp(\lambda \|b\|_{L^\infty}^2)\bigr)\leq \E\bigl(\exp(\lambda \|b\|_{C^\alpha}^2)\bigr)<\infty,
\ee
for some sufficiently small $\lambda>0$. This implies that $\E(\exp(k\|b\|_{L^\infty}))$ is
finite for all $0\leq k<\infty$. Note that when
$b$ has the representation \iref{paramb}, its covariance kernel is
then given by
\be
C_b(x,x')=\sum_{j\geq 1}\psi_j(x)\psi_j(x'),
\ee 
and therefore we may analyze the smoothness of $C_b$ through that of the individual $\psi_j$. 
A strategy to establish that $b$ with representation \iref{paramb} is almost surely in a 
H\"older space $C^\kappa$ for some $\kappa>0$ is discussed in \cite{H,DS}. It assumes that the $\psi_j$ are 
individually in $C^\alpha$ for some $\alpha>0$, and that for some $0<\delta<2$,
 \be
 \sum_{j\geq 1}\|\psi_j\|_{L^\infty}^2<\infty \quad {\rm and} \quad   \sum_{j\geq 1}\|\psi_j\|_{L^\infty}^{2-\delta}C_j^{\delta} <\infty,
 \label{cond2}
 \ee
where $C_j$ is a constant such that $|\psi_j(x)-\psi_j(x')|\leq C_j|x-x'|^\alpha$. Then, using the Kolmogorov
continuity theorem, one establishes that $b$ belongs almost surely to $C^\kappa$ for $\kappa<\alpha\delta/2$, 
see Corollary 7.22 in \cite{DS} (one can also easily check that \iref{cond2} implies \iref{cond1} with $\beta=\delta$). 

A second approach for guaranteeing Property \M{k}  does not assume any smoothness on the $\psi_j$, but instead the summability property
\be
\sum_{j\geq 1}\|\psi_j\|_{L^\infty}<\infty.
\label{cond3}
\ee
Then, as observed in \cite{HS}, we find that for all $0\leq k<\infty$,
\be
\E\bigl(\exp(k\|b(y)\|_{L^\infty})\bigr) \leq \E\(\exp\(k\sum_{j\geq 1} |y_j| \|\psi_j\|_{L^\infty}\)\) =\prod_{j\geq 1} \E\bigl(\exp(k|y_j| \|\psi_j\|_{L^\infty})\bigr) <\infty,
\ee
since for a standard Gaussian variable $t$ and a positive quantity $s$ one has 
\be
\E\bigl(\exp(s |t|)\bigr)=\frac 2{\sqrt {2\pi}} \int_{0}^{+\infty} \exp\biggl(s t-\frac {t^2}2\biggr)dt= e^{\frac {s^2}2}\frac 2{\sqrt {2\pi}}
\int_{-s}^{+\infty} e^{-\frac {t^2}2}dt\leq \exp\biggl(\frac {s^2}2+\frac 2{\sqrt {2\pi}}s\biggr).
\ee

Neither of these two types of conditions imply each other. Indeed, on the one hand, if $(\psi_j)_{j\geq 0}$ is a smooth 
wavelet basis of $L^2(D)$ ordered from coarser to finer scale, $C^\beta$ smoothness of $C_b$ holds 
whenever $\|\psi_j\|_{L^\infty} \leq C j^{-\frac \beta {2d}}$,
which is not sufficient to ensure that $(\|\psi_j\|_{L^\infty})_{j\geq 1}$ belongs to $\ell^1(\N)$ if $\frac \beta {2d}\leq 1$. 
This shows that \iref{cond1} may hold while \iref{cond3} may fail. Also note that, for the same reason,
the criterion \iref{cond2} is sufficient but not necessary for \iref{cond1} to hold.
On the other hand, if the $\psi_j$ are discontinuous, we cannot hope that $C_b$ has H\"older smoothness
while $(\|\psi_j\|_{L^\infty})_{j\geq 1}$ could belong to $\ell^1(\N)$, which shows that \iref{cond3} may hold
while \iref{cond1} and \iref{cond2} may fail. Note that the exclusion of discontinuous $\psi_j$ is problematic
when one needs to model sharp interfaces in the diffusion media.

We also remark that conditions \iref{cond2} and \iref{cond3} are artificially strong 
when the supports of the $\psi_j$ do not overlap, or only partially overlap.
This motivates us to  introduce a third condition which better takes into account the
support properties of the $\psi_j$, and does not enforce them to be H\"older continuous.
\nl
\nl
{\bf Assumption A:} {\it There exists a strictly positive sequence $\rho=(\rho_j)_{j\geq 1}$ such that 
the series $\sum_{j\geq 1} \rho_j |\psi_j|$ converges in $L^\infty(D)$ and  
\be\label{Asum}
\sum_{j\geq 1} \exp(-\rho_j^2)<\infty.
\ee
} 
\nl
Note that since the $\rho_j$ are strictly positive, the summability of the $\exp(-\rho_j^2)$ also implies that 
\be
\rho_0:=\inf_{j\geq 1}\rho_j >0.
\ee
We shall make use of the following elementary lemma.

\begin{lemma}
\label{lemtail}
If  $(\alpha_j)_{j\geq 1}$ is a sequence of   
numbers from $[0,1[$ that belongs to $\ell^1(\N)$, then,  there is a constant $c>0$ and $M_0\ge 1$ such that 
\be
\label{rem11}
\sum_{j=1}^\infty \alpha_j^M\le  e^{- c M},\quad  M\ge M_0.
\ee
\end{lemma}

\noindent
{\bf Proof:} We first consider the case when  $\sum_{j=1}^\infty \alpha_j<1$.  In this case, \iref{rem11} follows with $M_0=1$ because $\|(\alpha_j)_{j\geq 1}\|_{\ell^M}\le \| (\alpha_j)_{j\geq 1}\|_{\ell^1}$, for any $M\geq 1$.  Hence, with $c= -\ln \|(\alpha_j)_{j\ge 1}\|_{\ell^1}$, we arrive at \iref{rem11}.  To prove the general case, we first choose $K$ so that $\sum_{j> K}\alpha_j\le 1/2$.  If $M_0>1$ is sufficiently large then
 \be
 \label{rem111}
 \sum_{j=1}^\infty \alpha_j^{M_0} = \sum_{j=1}^K\alpha_j^{M_0} +\sum_{j> K}\alpha_j^{M_0}\le 3/4,
 \ee
 where we used the case already proven for estimating the second term.   Since $\alpha_j^M=(\alpha_j^{M_0})^{M/M_0}$ we know from the comparison of $\ell^p$ norms
 \be 
 \label{rem12}
  \sum_{j=1}^\infty \alpha_j^{M}\le  (3/4)^{M/M_0},\quad M\ge M_0.
 \ee
and the result follows. \hfill $\Box$
\nl

We now prove that Assumption {\rm \bf A} guarantees the measurability and finiteness of the exponential moments
of the map $y\mapsto b(y)$.  
\begin{theorem}\label{AimpliesMk}
Assumption {\rm \bf A} implies that Property \M{k}
holds for any $0\leq k<\infty$.
\end{theorem}

\noindent
{\bf Proof:}
Let us consider a sequence $(\rho_j)_{j\geq 1}$ such that \iref{Asum} holds. For any $t\geq 0$, the complement $\cE_t^c$ of the event
\be
\cE_t:=\bigl\{y \,: \, \sup_{j\geq 1} \rho_j^{-1} |y_j| \leq t \bigr\},
\ee
has measure
\be
\gamma(\cE_t^c)\leq  \sum_{j=1}^\infty\gamma\{y \, : \, |y_j|>t\rho_j\}\le 
 \frac {2}{\rho_0\sqrt {2\pi}} \sum_{j=1}^\infty  \exp\biggl(-\frac{t^2\rho_j^2}{2}\biggr), 
 \label{gammaE}
 \ee
 where we have used the univariate Gaussian tail bound
\be
\int_{|s|>B} g(s)ds \leq \frac {2}{B\sqrt {2\pi}} e^{-\frac {B^2} 2}, \quad  B>0.
\ee
Therefore, application of Lemma \ref{lemtail} with $\alpha_j=\exp(-\rho_j^2)$ shows that
\be
\gamma(\cE_t^c)\leq Ce^{-ct^2},  \quad t\geq t_0,
\ee
for certain constants $C,c,t_0$. In particular the event
\be
\cE:=\cup_{t\geq 0} \cE_t =\bigl\{y \,: \, \sup_{j\geq 1} \rho_j^{-1} |y_j| <\infty \bigr\},
\ee
has full measure, i.e.
\be
\gamma(\cE)=1.
\ee

We next introduce for $J\geq 1$ the truncation
\be
b_J(y)=\sum_{j=1}^J y_j\psi_j.
\label{truncb}
\ee
The mapping $y\mapsto b_J(y)$ is measurable from $U$ to $L^\infty(D)$ since it is
a continuous $L^\infty$-valued function of the variables $(y_1,\dots,y_J)$.  We observe that if Assumption {\bf A} holds, then, for any $y\in \cE$, we may
define $b(y)$ as the limit in $L^\infty(D)$ of $b_J(y)$ since
\be
 \|b_J(y)-b(y)\|_{L^\infty} \leq \Big \|\sum_{j>J} \rho_j |\psi_j| \Big \|_{L^\infty} \sup_{j\geq 1} \rho_j^{-1}|y_j|  \to 0 \quad {\rm as}\; J\to +\infty.
\ee
Thus $b_J(y)$ converges to $b(y)$ in $L^\infty(D)$ for almost every $y\in U$. Therefore $y\mapsto b(y)$ is 
also a measurable mapping, that is, $b(y)$ is a random variable with values in $L^\infty(D)$. 

In addition, 
for any $s\geq 0$ and $y\in \cE_s$, we may write
\be
\|b(y)\|_{L^\infty}\leq C_A s, \quad C_A:=\Big \|\sum_{j\geq 1} \rho_j |\psi_j| \Big \|_{L^\infty}.
\ee
Hence, for $t\ge 1$, we have
\be
\label{Athm1} 
P(t):=\gamma\{ y\, :\, \|b(y)\|_{L^\infty} >t\}\le 
\gamma(\cE_{t/C_A}^c) \leq
 \frac {2C_A}{\rho_0\sqrt {2\pi}} \sum_{j=1}^\infty  \exp\biggl(-\frac{t^2\rho_j^2}{2C_A^2}\biggr), 
\ee
where we have used \iref{gammaE}. From Assumption {\bf A}, we know that the last sum in \iref{Athm1} is finite when $t=t_0:=\sqrt{2}C_A$.
Hence, applying Lemma \ref{lemtail} with $\alpha_j= \exp\(-\frac{t_0^2\rho_j^2}{2C_A^2}\)$, and using the fact that $P(t)\leq 1$ for all $t\geq 0$, we find that 
 \be
 \label{Athm11} 
P(t)\le Ce^{-ct^2}, 
\ee
for suitable $c,C>0$ and for all $t\geq 0$. Therefore
\be
\E\bigl(\exp(k\|b(y)\|_{L^\infty})\bigr) =\int_{0}^{+\infty} k\exp(kt)P(t)dt <\infty,
\ee
for all $0\leq k<\infty$. \hfill $\Box$
\nl

If we next define $u_J(y)$ as the weak solution of \iref{pde} with $a=a_J(y)=\exp(b_J(y))$,
where $b_J$ is as in \iref{truncb},
we find that the mapping $y\mapsto u_J(y)$ is measurable from $U$ to $V$ since it is
a continuous $V$-valued function of the variables $(y_1,\dots,y_J)$.    This continuity  stems from the following
classical stability estimate: 
two weak solutions $u$ and $\t u$ of \iref{pde} with diffusion coefficients $a$ and $\t a$, respectively,
satisfy
\be
\|u-\t u\|_V \leq C \|a-\t a\|_{L^\infty}, \quad   C:=\frac {\|f\|_{V*}} {\min \{a_{\min},\t a_{\min}\}^2}.
\label{strang}
\ee
This estimate also shows
that the  convergence
of $b_J(y)$ towards $b(y)$ in $L^\infty(D)$ for each $y\in \cE$ implies the convergence of  $u_J(y)$ towards
$u(y)$ in $V$ for each $y\in \cE$, which shows that the mapping $y\mapsto u(y)$ is measurable as an almost everywhere
limit of measurable mappings,
 that is, $u(y)$ is a random variable with values in $V$.
 By application of \iref{laxmilgram}, we thus obtain the following result. 
 
\begin{corollary}\label{Aimpliesuk}
Assumption {\rm \bf A} implies that $y\mapsto u(y)$ is measurable with values in $V=H^1_0(D)$ 
and that $ \E\bigl(\|u(y)\|_V^k\bigr)$ is finite for
any $0\leq k<\infty$.
\end{corollary}

\begin{remark}
Assumption {\bf A} is almost necessary for Property \M{k} to hold 
in the case where the supports of $\psi_j$ do not overlap, in the sense that Property  
\M{k} then implies that $\sum_{j\geq 1}\rho_j|\psi_j|$ is uniformly bounded.
In this case, the uniform boundedness
of $\sum_{j\geq 1}\rho_j|\psi_j|$ for a sequence $(\rho_j)_{j\geq 1}$ which satisfies \iref{Asum}
can be equivalently expressed
as follows: there
exists $C_B>0$ such that with $b_j:=\|\psi_j\|_{L^\infty}$, one has
\be
\sum_{j\geq 1} \exp\biggl(-\frac {C_B}{b_j^2}\biggr) <\infty.
\label{equiv}
\ee
Indeed, since $C_A:=\Big \|\sum_{j\geq 1}\rho_j |\psi_j| \Big \|_{L^\infty}=\sup_{j\geq 1} \rho_j b_j$,  we have \iref{equiv} with $C_B=C_A^{-2}$.
We then have
\be
\gamma\{ y\, :\, \|b(y)\|_{L^\infty(D)} \leq t\}= \prod_{j\geq 1}\gamma\Bigl\{ y\; :\; |y_j| \leq \frac {t}{b_j}\Bigr\} \leq \prod_{j\geq 1}\biggl(1-\frac {2t}{b_j\sqrt {2\pi}}\exp\biggl(-2 \frac {t^2}{b_j^2}\biggr)\biggr),
\ee
for any $t>0$, where we have used the lower bound
\be
\int_{|s|>B} g(s)\,ds \geq \int_{B< |s|< 2B} g(s)\,ds \geq\frac {2B}{\sqrt {2\pi}}\exp(-2 B^2).
\ee
It follows that if \iref{equiv} does not hold for any $C_B>0$, then $\gamma\{ y\, :\, \|b(y)\|_{L^\infty(D)} \leq t\}=0$ for all $t>0$, which
means that $\|b(y)\|_{L^\infty(D)}=\infty$ with probability $1$ and therefore Property \M{k} does not hold.
\end{remark}

\begin{remark}\label{Aisimplied}
Assumption {\bf A} always follows from the 
assumptions of Theorem \ref{maintheo}. Indeed, if $(\rho_j^{-1})_{j\geq 1}\in \ell^q(\N)$
for some $q<\infty$, one can find a positive sequence $(\delta_j)_{j\geq 1}$ such that
$\lim_{j\to +\infty}\delta_j=0$ and such that the sequence defined as $\omega_j:=\delta_j \rho_j$ satisfies
$\bigl(\exp(-\omega_j^2)\bigr)_{j\geq 1} \in \ell^1(\N)$. We then obtain the validity of Assumption {\bf A}
with the sequence $(\omega_j)_{j\geq 0}$ since
\be
\Big\| \sum_{j>J} \omega_j|\psi_j| \Big\|_{L^\infty} \leq \(\sup_{j>J} \delta_j \)\Big\| \sum_{j\geq 1} \rho_j|\psi_j| \Big \|_{L^\infty} \to 0, \quad {\rm as}\quad J\to +\infty.
\ee
\end{remark}

\section{Hermite expansions and partial derivatives}

We are interested in the summability properties of the multivariate Hermite expansion
\iref{multiherm}. By Parseval's identity, the $L^2$ integrability of $u(y)$ with respect to
the Gaussian measure $\gamma$ implies the $\ell^2$ summability property
\be
\sum_{\nu\in \cF} \|u_\nu\|_V^2 <\infty.
\ee
In order to obtain polynomial approximation results in the mean-square sense, or equivalently, in $L^2(\R^\N,V,\gamma)$, we 
need to establish the $\ell^p$ summability 
\be
\sum_{\nu\in \cF} \|u_\nu\|_V^p <\infty,
\ee
for values of $p<2$. For this purpose, we  establish   $L^2$-integrability properties of the 
partial derivatives of $u$.  For any $V$-valued function $y\mapsto w(y)$ defined on $\R^\N$, and $j\in\N$, the derivative $\partial_{y_j}w$ is defined as the limit, as $h\to 0$,  of the difference quotient
\be
\label{wder}
 \frac{ w(y+he_j)-w(y)}{h}
\ee
provided this limit exists in $V$, where
\be
e_j:=(\delta_{i,j})_{i\geq 1}
\ee
is the Kronecker sequence of index $j$. Higher  derivatives $\partial^\nu w$ are defined inductively.  For $\nu\in \cF$, we use the notation
\be
\partial^\nu w=\(\prod_{j\geq 1} \partial_{y_j}^{\nu_j}\)w .
\ee
      When, $\nu_j=0$, the operator $ \partial_{y_j}^{\nu_j}$ is the identity operator and so the above product
has only a finite number of factors which are not the identity.   

 We also use the standard notation ${\rm supp}(\nu):=\{j\; : \; \nu_j\neq 0\}$, as well as
\be  
|\nu|:=\|\nu\|_{\ell^1}=\sum_{j\geq 1} \nu_j,\quad \nu\in \cF,
\ee
and
\be
\nu !:=\prod_{j\geq 1}\nu_j !,\quad  \nu\in\cF,
\ee
with the convention that $0!=1$.   Furthermore, we use the combinatorial notation,
\be
{\nu\choose \mu}:=\prod_{j\geq 1}{\nu_j \choose \mu_j},\quad \mu,\nu\in\cF,
\ee
with the convention that
\be
\label{choosezero}
{n\choose m}:=0, \quad {\rm if}\quad m>n.
\ee

In  analogy with $y=(y_j)_{j\geq 1}$, we set $\psi=(\psi_j)_{j\geq 1}$ and use the notation
\be
\psi^\nu:=\prod_{j\geq 1} \psi_j^{\nu_j}, \quad \nu\in\cF.
\ee 
Finally, with $\leq$ denoting the componentwise partial order between multi-indices, we define
\be
\label{shadow}
S_\mu:=\{\nu \in \cF: \nu \le \mu \ {\rm and } \ \nu \neq \mu\},\quad \mu\in\cF,
\ee
which is always a finite set.
The following lemma proves  that the derivatives $\partial^\nu u(y)$ are elements of $V$ and gives a recursive way of computing them.
While it is implicitly contained in the proof of Theorem 3.1 in \cite{HS}, 
we give it here in an explicit form for completeness.

\begin{lemma}
For any $y$ such that $\|b(y)\|_{L^\infty}<\infty$, and any nonzero $\nu\in\cF$, the 
partial derivative $\partial^\nu u(y)$ exists and is the solution of the variational problem
\be
\int_D a(y)\nabla \partial^\mu u(y)\cdot \nabla v =-\sum_{\nu \in S_\mu} {\mu\choose \nu}\int_D \psi^{\mu-\nu} a(y)  \nabla \partial^\nu u(y)\cdot \nabla v , \quad v\in V.
\label {recurs}
\ee
\end{lemma}

\noindent
{\bf Proof:} The proof of this lemma is, in principle, the same as the  proof of Theorem 4.2 in \cite{CDS}.  The only difference is that
in the theorem of \cite{CDS} the assumption is that $a$ is affine in $y$,  whereas now we assume $a=e^b$ with $b$ affine
in $y$.   In going further, we  indicate only the necessary changes  caused by the different diffusion coefficient.  

One begins by establishing the validity of the theorem when $\mu=e_j$, for any $j\in\N$, by considering the functions $\rho_h(y):=\frac{u(y+he_j)-u(y)}{h}$, for $|h|\le 1$.
Following the reasoning in \cite{CDS},  one derives that $\rho_h$ satisfies
\be
\label{satisfies}
\int_D a(y)\nabla \rho_h(y)\cdot \nabla v = L_h(v),\quad 
\mbox{for all } v\in V,
\ee
where $L_h$ is the linear functional on $V$ given by 
\be
L_h(v):=-\int_D \frac{e^{h\psi_j}-1}{h}a(y)\nabla u(y+he_j)\cdot \nabla v.
\ee
Using \iref{strang}, one next proves that $L_h\to L$ in $V'$ where $L$ is the linear functional
\be
L(v):= -\int_D  a(y)\psi_j \nabla u(y)\cdot \nabla v.
\ee
Therefore, taking a limit as $h\to 0$ in \iref{satisfies}  proves the theorem for the case $\mu=e_j$.  The case of 
any general $\mu\in\cF$ now follows by a recursive application of this  result.	\hfill $\Box$
\nl

The recursive formula \iref{recurs} will be crucial for obtaining sharp estimates
for the partial derivatives $\partial^\mu u(y)$.  For the time 
being, we use it to establish the  following crude bound
\be
\| \partial^\mu u(y) \|_V \leq C\exp\bigl(( 2 | \mu|+1) \|b(y)\|_{L^\infty}\bigr),  \quad \mu \in \cF,
\label{crude}
\ee
where the constant depends on $\mu$, on the $\psi_j$, and on $f$. In view of  Property \M{k}, this guarantees the  finiteness of all
moments $\E(\| \partial^\mu u(y) \|_V^k)$ for all $0\leq k<\infty$.

For the proof of \eqref{crude}, for each $y$ such that $\|b(y)\|_{L^\infty}<\infty$, we introduce a $y$-dependent equivalent norm on $V$,
\be
   \| v \|_{a(y)}^2 := \int_D a(y) | \nabla v |^2 \, dx ,\quad v \in V\,. 
   \label{ay}
\ee
For such $y$, we have the norm equivalence
\be\label{normequiv}
\exp\bigl(-\|b(y)\|_{L^\infty}\bigr)\|v\|_{V}^2 \leq \|v\|_{a(y)}^2 \leq  \exp\bigl(\|b(y)\|_{L^\infty}\bigr) \|v\|_{V}^2.
\ee
Using this in conjunction with  \eqref{recurs} and $v=\partial^\mu u(y)$, we obtain
\begin{align*}
   \| \partial^\mu u(y)\|_V^2 &\leq e^{\|b(y)\|_{L^\infty}}  \| \partial^\mu u(y) \|^2_{a(y)}   \\
   & \leq e^{\|b(y)\|_{L^\infty}} \sum_{\nu \in S_\mu}  {\mu\choose \nu} \(\prod_{j\in \operatorname{supp} (\mu-\nu)} \| \psi_j\|_{L^\infty}^{\mu_j  - \nu_j}\)   \| \partial^\nu u(y)\|_{a(y)}  \| \partial^\mu u(y)\|_{a(y)}   \\
    & \leq  C_\mu  \,  e^{2 \|b(y)\|_{L^\infty}}   \| \partial^\mu u(y)\|_V  \sum_{\nu \in S_\mu} \| \partial^\nu u(y)\|_V,
 \end{align*}
	 with a constant $C_\mu>0$.  
	 Applying this recursively, starting from the estimate \eqref{laxmilgram} for $\| u(y)\|_V$, we obtain \eqref{crude}.

The following result relates the $L^2$ norms of the
mixed derivatives of $u$ up to some given order $r$ with 
weighted $\ell^2$ norms of the Hermite coefficients.

\begin{theorem}
\label{identity}
Let $r\geq 0$ be an integer and let $\rho=(\rho_j)_{j\geq 1}$ be a sequence of positive numbers. Assume that Property \M{k} hold for all $0\leq k<\infty$. Then
\be\label{identity1}
\sum_{\|\mu\|_{\ell^\infty} \leq r}\frac {\rho^{2\mu}}{\mu !} \int_U \|\partial^\mu u(y)\|_V^2 \, d\gamma(y)=\sum_{\nu\in \cF} b_\nu \|u_\nu\|_V^2,
\ee
where
\be\label{identity2}
b_\nu:=\sum_{\|\mu\|_{\ell^\infty}\leq r}{\nu\choose \mu} \rho^{2\mu}.
\ee
The weights $b_\nu$ can also be expressed as
\be\label{identity3}
b_\nu=\prod_{j\geq 1}\(\sum_{l=0}^r {\nu_j\choose l}\rho_j^{2l}\).
\ee
\end{theorem}

\noindent
{\bf Proof:}
With the normalization in $L^2(\R, dg)$ that we have chosen for the univariate Hermite polynomials, we have the Rodrigues formula
\be
   H_n(t) =  \frac{(-1)^n}{\sqrt{n!} } \frac{g^{(n)}(t)}{g(t)} \,.
\ee
As a consequence, for $m \le n$ and sufficiently smooth univariate functions $v$ with square-integrable derivatives, integration by parts yields
\be\label{ibp1d}
  v_n :=\int_\R v \,H_n \,g \,dt = \frac{(-1)^n}{\sqrt{n!}} \int_\R v \, g^{(n)} \,dt    = \sqrt{\frac{(n-m)!}{n!}}   \int_\R v^{(m)} \, H_{n-m} \, g\,dt.
  \ee
By Parseval's identity, we thus have
\be
  \frac{1}{m!} \int_\R | v^{(m)}|^2 \, g\,dt = \sum_{n \geq m} \frac{n!}{m! (n-m)!} |v_n|^2 = \sum_{n\geq 0} {n\choose m} |v_n|^2 ,
\ee
where in the last step we have made use of our notational convention \eqref{choosezero}.

For $\nu,\mu\in\cF$ with $\mu \leq \nu$, applying the same reasoning to the tensor product Hermite coefficients $u_\nu$, we obtain
\be
 u_\nu = \langle u, H_\nu\rangle = \sqrt{ \frac{(\nu - \mu)!} {\nu!} } \, \langle \partial^\mu u, H_{\nu - \mu} \rangle \,.
\ee
Concerning the integration by parts, note that, as a consequence of \iref{crude}, Property \M{k} allows us to proceed as in \eqref{ibp1d}.
By Parseval's identity, summation over $\nu \in \cF$ with $\nu \geq \mu$ thus gives
\be
   \frac1{\mu!} \int_U \| \partial^\mu u(y)\|_V^2 \,d\gamma(y) = \sum_{\nu \geq \mu}  \frac{\nu!}{\mu!(\nu-\mu)!} \| u_\nu \|^2_V = \sum_{\nu \in \cF} {\nu\choose\mu} \| u_\nu \|^2_V, 
\ee
where we have again used \eqref{choosezero}.
Multiplying this identity by $\rho^{2\mu}$ and summing over $\mu$ with $\|\mu\|_{\ell^\infty}\le r$ gives
\be
   \sum_{\| \mu \|_{\ell^\infty} \le r}  \frac{\rho^{2\mu}}{\mu!} \int_U \| \partial^\mu u(y)\|_V^2 \,d\gamma(y) =    \sum_{\| \mu \|_{\ell^\infty} \le r}   \sum_{\nu \in \cF} {\nu\choose\mu} \rho^{2\mu} \| u_\nu \|^2_V = \sum_{\nu \in \cF} b_\nu \| u_\nu \|^2_V,
\ee
that is, \eqref{identity1} with $b_\nu$ as defined in \eqref{identity2}.
The equivalent expression of $b_\nu$,  given in  \eqref{identity3},
immediately follows by factorization in \eqref{identity2}.
\hfill $\Box$
\nl

\section{Estimates of the partial derivatives}

In view of Theorem \ref{identity}, estimating certain
weighted $\ell^2$ norms of the sequence
$(\|u_\nu\|_V)_{\nu\in \cF}$ amounts to estimating related weighted $\ell^2$ norms
of the sequence
\be
\biggl( \int_U \|\partial^\mu u(y)\|_V^2 \, d\gamma(y)\biggr)_{\|\mu\|_{\ell^\infty} \leq r},
\ee
which we next address.

\subsection{Estimates for fixed $y$}

Recalling the norm $\|\cdot\|_{a(y)}$ defined in \iref{ay}, we first fix $y$ and estimate the relevant weighted $\ell^2$ norms
for the sequence $\bigl(\|\partial^\mu u(y)\|_{a(y)}\bigr)_{\|\mu\|_{\ell^\infty} \leq r}$.

\begin{theorem}
\label{theofixedy}
Let $r\geq 1$ be an integer. Assume that there exists a positive sequence $(\rho_j)_{j\geq 1}$ such that
\be
\label{weightcondition}
\sup_{x\in D} \sum_{j\geq 1}\rho_j|\psi_j(x)| =:K < C_r:=\frac {\ln 2}{\sqrt r}.
\ee
Then, there exists a constant $C$ that depends on $K$ and $r$, such that
\be
\sum_{\|\mu\|_\infty \leq r} \frac {\rho^{2\mu}}{\mu !} \|\partial^\mu u(y)\|_{a(y)}^2  \leq C \| u(y)\|_{a(y)}^2,
\ee
for all $y$ satisfying $\|b(y)\|_{L^\infty}<\infty$.
\end{theorem}

\noindent
{\bf Proof:} For $k\geq 0$, we define  $\Lambda_k:=\{\mu\in \cF : |\mu|=k, \ \|\mu\|_{\ell_\infty}\le r\}$ and 
\be
\label{defsigmak}
\sigma_k :=\sum_{\mu\in\Lambda_k}\frac {\rho^{2\mu}} {\mu ! } \|\partial^\mu u(y)\|_{a(y)}^2. 
\ee 
We prove that 
\be
\label{lambdak}
\sigma_k  
\le \sigma_0\delta^k,
\ee
for a fixed $\delta<1$.  Since $\sigma_0=\| u(y)\|_{a(y)}^2$,  the theorem follows from this  by summing over $k$.

 We introduce the notation
\be
\label{notation}
\e(\mu,\nu):=\frac{\sqrt{\mu !}}{\sqrt{\nu !}}\frac {\rho^{\mu-\nu}|\psi| ^{\mu-\nu}}{(\mu-\nu)!},
\ee
where $|\psi|^\nu:=\prod_{j\geq 1} |\psi_j|^{\nu_j}$. From \iref{recurs}, recalling the notation \eqref{shadow}, we have  
\be
\sigma_k\le   \int_D\sum_{\mu\in\Lambda_k}   \sum_{\nu\in S_\mu} \e(\mu,\nu)a(y) \frac{\rho^\nu |\nabla \partial^{\nu}u(y)|}{\sqrt{\nu !}} \frac{\rho^\mu|\nabla \partial^{\mu}u(y)|}  { \sqrt{\mu !}} ,
\ee
and thus, by Cauchy-Schwarz inequality,
\be
 \label{Dnu21}
\sigma_k\le   \int_D   \sum_{\mu\in\Lambda_k} \biggl(\sum_{\nu\in S_\mu}     \e(\mu,\nu)a(y) \frac{|\rho^\nu\nabla \partial^{\nu}u(y)|^2}{\nu !} \biggr)^{1/2}
 \biggl( \sum_{\nu\in S_\mu}      \e(\mu,\nu) a(y) \frac{|\rho^\mu\nabla \partial^{\mu}u(y)|^2}{ \mu !} \biggr)^{1/2}.
\ee
 For $\mu\in \Lambda_k$,
 define $S_{\mu,\ell}:= \{ \nu\in S_\mu : \ |\mu-\nu|=\ell\}$. We have
\be
\label{crucial1}
 \sum_{\nu\in S_\mu}\e(\mu,\nu)=\sum_{\ell =1}^k\sum_ {\nu\in S_{\mu,\ell}}\e(\mu,\nu).
\ee
Now,  from \iref{weightcondition}, we have
\begin{align*}
\sum_ {\nu\in S_{\mu,\ell}}\e(\mu,\nu)&\le  r^{\ell/2}\sum_{\nu\in S_{\mu,\ell}} \frac{ \rho^{\mu-\nu}|\psi| ^{\mu-\nu}}{(\mu-\nu)!}
\leq r^{\ell/2}\sum_{|\tau|=\ell} \frac{ \rho^{\tau}|\psi| ^{\tau}}{\tau!}\\
&=r^{\ell/2}\(\sum_{j=1}^\infty \rho_j |\psi_j(x)|\)^\ell= r^{\ell/2}\frac{1}{\ell!}K^\ell.
\end{align*}
Inserting this into \iref{crucial1} we find that 
\be
\label{find}
\sum_{\nu\in S_\mu}\e(\mu,\nu) \le e^{\sqrt{r}K}-1\le e^{\sqrt{r}C_r}-1 = 1.
\ee
Inserting this into  \iref{Dnu21}, and applying the Cauchy-Schwarz inequality, we obtain
\begin{align*}
\sigma_k
&\le   \int_D   \sum_{\mu\in\Lambda_k} \biggl( \sum_{\nu\in S_\mu}     \e(\mu,\nu)a(y) \frac{|\rho^\nu\nabla \partial^{\nu}u(y)|^2}{\nu !} \biggr)^{1/2} 
\biggl( a(y) \frac{|\rho^\mu\nabla \partial^{\mu}u(y)|^2}{\mu !} \biggr)^{1/2} \\
&\le   \int_D \biggl(\sum_{\mu\in\Lambda_k} \sum_{\nu\in S_\mu}     \e(\mu,\nu)a(y) \frac{|\rho^\nu\nabla \partial^{\nu}u(y)|^2}{\nu !} \biggr)^{1/2} 
 \biggl(\sum_{\mu\in\Lambda_k} a(y) \frac{|\rho^\mu\nabla \partial^{\mu}u(y)|^2}{\mu !} \biggr)^{1/2}.
\end{align*}
We treat the first factor in the last integral by interchanging summations in $\mu$ and $\nu$.
For this purpose, we introduce, for any $\ell \leq k-1$ and $\nu\in \Lambda_{\ell}$, the set
\be
R_{\nu,k}:=\{\mu\in \Lambda_k : \nu\in S_\mu\}.
\ee
By a similar argument as above for the set $S_{\mu,\ell}$, we find that 
\be
\sum_{\mu\in R_{\nu,k}} \e(\mu,\nu)\leq r^{(k-\ell)/2}\frac{1}{(k-\ell)!}K^{k-\ell}.
\ee
It follows that
\be
\sigma_k\leq  \int_D   \biggl(  \sum_{\ell=0}^{k-1}r^{(k-\ell)/2}\frac{1}{(k-\ell)!}K^{k-\ell} \sum_{\nu\in \Lambda_\ell}    a(y) \frac{|\rho^\nu\nabla \partial^{\nu}u(y)|^2}{\nu !} \biggr)^{1/2} 
\biggl(\sum_{\mu\in\Lambda_k} a(y) \frac{|\rho^\mu\nabla \partial^{\mu}u(y)|^2}{\mu !} \biggr)^{1/2} .
\ee
If we now apply the Cauchy-Schwarz inequality on the integral we obtain
 \be  \label{Dnu4}
\sigma_k\le  \( \sum_{\ell=0}^{k-1}\frac{1}{(k-\ell)!}(\sqrt{r}K)^{k-\ell} \sigma_\ell\)^{1/2}\sigma_k^{1/2}.
 \ee
In other words,
\be  \label{Dnu5}
\sigma_k\le \sum_{\ell=0}^{k-1}\frac{1}{(k-\ell)!}(\sqrt{r}K)^{k-\ell} \sigma_\ell .
 \ee
 Now pick $\delta<1$ such that $K/\delta\leq C_r$.  By induction,  we prove that $\sigma_k\leq \sigma_0\delta^k$ for all $k\geq 0$.   This is clearly true for $k=0$ and assuming it has been proven for $k-1$, we find   
\be
\sigma_k\le \sigma_0\sum_{\ell=0}^{k-1}\frac{1}{(k-\ell)!}(\sqrt{r}K)^{k-\ell} \delta^\ell \leq \sigma_0\delta^k \sum_{\ell=0}^{k-1}\frac{1}{(k-\ell)!}(\sqrt{r}C_r)^{k-\ell}
\le \sigma_0 \delta^k( e^{\sqrt{r}C_r}-1 )= \sigma_0 \delta^k,
\ee
which advances the induction.
This establishes \iref{lambdak} and completes the proof of the theorem. \hfill $\Box$
\vskip .2in

\subsection{Integral estimates}

Combining \eqref{normequiv} with Theorem \ref{theofixedy}, we obtain the following result.

\begin{theorem}
\label{theoint}
Let $r\geq 1$ be an integer. Assume that there exists a positive sequence $(\rho_j)_{j\geq 1}$ such that
\be
\sup_{x\in D} \sum_{j\geq 1}\rho_j|\psi_j(x)| =:K < C_r:=\frac {\ln 2}{\sqrt r}.
\ee
Assume in addition that Property \M{4} holds.
Then we have
\be
\sum_{\|\mu\|_{\ell^\infty} \leq r}\frac {\rho^{2\mu}}{\mu !} \int_U \|\partial^\mu u(y)\|_V^2 
\,d\gamma(y)<\infty .
\ee
\end{theorem}

\noindent
{\bf Proof:}  We apply \iref{normequiv} with $v= \partial^\mu u(y)$ and obtain
\be\label{ne}
\|\partial^\mu  u(y)\|_{V}^2 \leq \exp\bigl(\|b(y)\|_{L^\infty}\bigr) \|\partial^\mu u(y)\|_{a(y)}^2.
\ee
This gives
\begin{eqnarray}
\label{gives1}
\sum_{\|\mu\|_{\ell^\infty} \leq r}\frac {\rho^{2\mu}}{\mu !} \int_U \|\partial^\mu u(y)\|_V^2 \, d\gamma(y)&\le&  \int_U  \exp\bigl(\|b(y)\|_{L^\infty}\bigr)\sum_{\|\mu\|_{\ell^\infty} \leq r}\frac {\rho^{2\mu}}{\mu !} \|\partial^\mu u(y)\|_{a(y)}^2 \, d\gamma(y)\nonumber \\
&\le& C_1 \int_U   \exp\bigl(\|b(y)\|_{L^\infty}\bigr) \|u(y)\|^2_{a(y)} \, d\gamma(y)\nonumber\\
&\le& C_1 \int_U   \exp\bigl(2\|b(y)\|_{L^\infty}\bigr) \|u(y)\|^2_{V} \, d\gamma(y)\nonumber\\
&\le&  C_2 \int_U  \exp\bigl(4\|b(y)\|_{L^\infty}\bigr)   \, d\gamma(y),
\end{eqnarray}
where the second inequality uses Theorem \ref{theofixedy}, the third  inequality uses  \iref{normequiv} again, and the last inequality uses
\iref{laxmilgram}.  From Property \M{4}, the last integral is finite and the proof of the theorem is complete.
\hfill $\Box$
\nl

\section{Summability of Hermite coefficients}

According to Theorem \ref{identity}, the conclusion of Theorem \ref{theoint} also gives the weighted $\ell^2$ summability estimate
\be
\sum_{\nu\in \cF} b_\nu \|u_\nu\|_V^2 <\infty.
\label{weightedell2}
\ee
Using H\"older's inquality, for $0<p<2$ and $q:=\frac {2p}{2-p}$ we thus have 
\be\label{hoelder}
\sum_{\nu\in \cF}\|u_\nu\|_V^p\leq \(\sum_{\nu\in \cF} b_\nu \|u_\nu\|_V^2\)^{p/2} \(\sum_{\nu\in \cF} b_\nu^{-q/2}\)^{1-p/2}.
\ee
The following result allows us to control the second factor.

\begin{lemma}
\label{bsum}
Let $0<p<2$ and  $q:=\frac {2p}{2-p}$ and let  $(\rho_j)_{j\geq 1}$ be a positive sequence such that
\be
(\rho_j^{-1})_{j\geq 1} \in \ell^q(\N).
\ee
Then, for any positive integer $r$ such that $\frac 2 {r+1}<p$, the sequence $(b_\nu)$ 
defined in \eqref{identity2} for this $r$ satisfies
\be
\sum_{\nu\in \cF} b_\nu^{-q/2}<\infty.
\ee
\end{lemma}

\noindent
{\bf Proof:}
The sum under consideration can be rewritten in factorized form as
\be
\label{identity4}
  \sum_{\nu \in \cF} b_\nu^{-q/2} =  \sum_{\nu \in \cF} \prod_{j\geq 1} \( \sum_{l=0}^r {\nu_j\choose l}\rho_j^{2l}   \)^{-q/2}
    = \prod_{j\geq 1} \sum_{n \geq 0} \( \sum_{l=0}^r {n \choose l}\rho_j^{2l}   \)^{-q/2},
\ee
provided we can show that the product on the right side is finite. Now,
\be\label{bsumest}
  \sum_{n \geq 0} \( \sum_{l=0}^r {n \choose l}\rho_j^{2l}   \)^{-q/2} 
    \leq  \sum_{n \geq 0 }  \biggl[  {n \choose {n\wedge r}} \rho_j^{2(n\wedge r)}  \biggr]^{-q/2}
      \leq 1 + \rho_j^{-q} + \ldots + \rho_j^{-(r-1)q} + C_{r,q} \rho_j^{-rq},
\ee
with the constant
\be
 C_{r,q} := \sum_{n \geq r} {n\choose r}^{-q/2} = (r!)^{q/2} \sum_{n\geq 0} \bigl[(n+1)\cdots (n+r)\bigr]^{-q/2}.
\ee
We have $C_{r,q}< \infty$ if and only if $q> 2/r$, which holds precisely under our assumption $p>\frac2{r+1}$.
Since $\rho_j\to\infty$, there exists $J$ such that $\rho_j > 1$ for all $j> J$. For such $j$, we can bound the right side of \eqref{bsumest}
by $1 + (C_{r,q} + r - 1) \rho_j^{-q}$. Hence, returning  to \eqref{identity4}, for a finite constant $C$ depending on $w$, we have  
 \be
   \sum_{\nu \in \cF} b_\nu^{-q/2} \leq C\prod_{j > J} \bigr( 1 + (C_{r,q} + r - 1) \rho_j^{-q}\bigr),
 \ee
 where the product on the right side converges because $\sum_{j \geq 1} \rho_j^{-q}<\infty$.
\hfill $\Box$
\nl

We are now ready to prove our main result.
\nl

\noindent
{\bf Proof of Theorem \ref{maintheo}:} By our assumptions, we are given a positive sequence $(\rho_j)_{j\geq 1}$ such that
\be
\sup_{x\in D}\sum_{j\geq 1} \rho_j|\psi_j(x)| <\infty \quad {\rm and}\quad
(\rho_j^{-1})_{j\geq 1} \in \ell^q(\N),\quad q:=\frac {2p}{2-p}.
\label{w}
\ee
We choose $r$ as the minimal integer satisfying $\frac{2}{r+1}< p$.
Since the sequence $\rho_j$ can be rescaled by multiplying by an arbitrary constant without affecting \iref{w}, we can assume without loss of generality that
\be
  \sup_{x\in D} \sum_{j\geq 1}\rho_j|\psi_j(x)|  <  \frac {\ln 2}{\sqrt r}.
\ee
Since, as noted in Remark \ref{Aisimplied}, $(\rho_j^{-1})_{j\geq 1} \in \ell^q(\N)$ implies 
Assumption {\bf A}, by Theorem \ref{AimpliesMk} we consequently also have Property \M{k} for all $k$. We can thus combine Theorems \ref{identity} and \ref{theoint} to obtain
\be
 \sum_{\nu\in \cF} b_\nu \|u_\nu\|_V^2   =  \sum_{\|\mu\|_{\ell^\infty} \leq r}\frac {\rho^{2\mu}}{\mu !} \int_U \|\partial^\mu u(y)\|_V^2 \,d\gamma(y)<\infty .
\ee
Using this together with Lemma \ref{bsum} in \eqref{hoelder}, we arrive at the conclusion.
\hfill $\Box$

\begin{remark}
\label{remfin}
As stated in Theorem \ref{maintheo}, application
of Stechkin's lemma gives that the 
best $n$-term approximation polynomials
\be
u_{\Lambda_n}:=\sum_{\nu\in \Lambda_n} u_\nu H_\nu,
\ee
obtained by retaining the indices of the $n$ largest $\|u_\nu\|_V$,
satisfy the estimate 
\be
\|u-u_{\Lambda_n}\|_{L^2(U,V,\gamma)}\lsim n^{-s},
\ee
where $s:=\frac 1 p-\frac 1 2=\frac 1 q$. There is, however, a more direct
and constructive way of retrieving this convergence rate,
namely taking instead $\Lambda_n$ to be the set of indices corresponding to the $n$ smallest
values of the weights $b_\nu$ which appear in \iref{weightedell2}.
We then directly obtain that
\be
\|u-u_{\Lambda_n}\|_{L^2(U,V,\gamma)}
\leq   \sup_{\nu\notin\Lambda_n} b_\nu^{-1/2} \(\sum_{\nu\in\cF} b_\nu \|u_\nu\|_V^2\)^{1/2}
\lsim d_{n+1}^*,
\ee
where $(d_n^*)_{n\geq 1}$ is the decreasing rearrangement of the sequence $(b_\nu^{-1/2})_{\nu\in \cF}$.
As seen in Lemma \ref{bsum}, this sequence
belongs to $\ell^q(\cF)$ which implies that $d_n^*\lsim n^{-s}$ with $s:=\frac 1 q$.
\end{remark}

\section{Examples}

In this section, we present several examples of applications of Theorem \ref{maintheo} corresponding to
different support properties of the $(\psi_j)_{j\geq 1}$. In each case, we discuss which
range of $\ell^q$ summability 
of the sequence $(\|\psi_j\|_{j\geq 1})_{j\geq 1}$ implies $\ell^p$ summability of the 
sequence $(\|u_\nu\|_V)_{\nu\in\cF}$ for some $p<2$.

\subsection{Finitely overlapping supports}

We say that the family $(\psi_j)_{j\geq 1}$ has finitely overlapping supports if and only if there
exists an integer $M$ such that for every $x\in D$, 
\be
\#\{j \, : \, \psi_j(x)\neq 0\} \leq M.
\ee
One example with $M=1$, which corresponds to disjoint supports, is the set of characteristic functions 
\be
\psi_j=c_j\Chi_{D_j},
\ee
with some normalizing factor $c_j$, when $(D_j)_{j\geq 1}$ is a partition of $D$. 
Another example with $M\geq 1$ is the set of Lagrange finite element basis functions
of a given order $k\geq 1$, associated to a conforming simplicial partition of $D$.

For such families, we find that the choice
\be
\rho_j^{-1}:=\|\psi_j\|_{L^\infty},
\ee
yields
\be
\sup_{x\in D}\sum_{j\geq 1} \rho_j|\psi_j(x)| \leq M,
\ee
and therefore condition \iref{firstcond} in Theorem \ref{maintheo} is satisfied. We thus obtain
the following immediate corollary.

\begin{corollary}
Let $(\psi_j)_{j\geq 1}$ be a family with finitely overlapping supports, and let $0<p<2$. 
If $(\|\psi_j\|_{L^\infty})_{j\geq 1}$ belongs to $\ell^q(\N)$ for $q=q(p):=\frac {2p}{2-p}$,
then $(\|u_\nu\|_V)_{\nu\in \cF}$ belongs to $\ell^p(\cF)$. In particular, best $n$-term Hermite approximations
converge in $L^2(U,V,\gamma)$ with rate $n^{-s}$ where $s=\frac 1 p-\frac 1 2=\frac 1 q$.
\end{corollary}

As already observed, we always have $q(p)>p$, which shows that there is in this case
an improvement in  the summability properties of $(\|u_\nu\|_V)_{\nu\in \cF}$ 
over those of  $(\|\psi_j\|_{L^\infty})_{j\geq 1}$. For example, $\ell^2$ summability of $(\|\psi_j\|_{L^\infty})_{j\geq 1}$
 implies $\ell^1$ summability of $(\|u_\nu\|_V)_{\nu\in \cF}$ , and therefore
convergence of best $n$-term Hermite approximations with rate $n^{-1/2}$.

\subsection{Wavelets}

For a general wavelet bases on a domain $D\subset \R^d$,
we adopt the notation $(\psi_\lambda)$, used for example in \cite{Co}, where $\lambda$ concatenates the 
scale and spatial indices, with the convention that  the scale level $l$
of $\psi_\lambda$ is denoted by $|\lambda|$, i.e., $|\lambda|:= l$. Thus, there
are $\cO(2^{dl})$ wavelets at level $l$ and 
the wavelets at each given scale have finite overlap.  This means that 
for all $x\in D$,
\be
\#\{\lambda \, : \, |\lambda|=l\; {\rm and}\; \psi_\lambda(x)\neq 0\} \leq M,
\ee
for some fixed $M>0$ independent of $l$. We consider wavelets
normalized such that
\be
\|\psi_\lambda\|_{L^\infty} =  c_l=C2^{-\alpha l}, \quad |\lambda|=l,
\label{normalpha}
\ee
for some fixed $C>0$ and $\alpha>0$. Using the finite overlapping property, 
we find that for any $0<\kappa<\alpha$ the sequence 
\be
\label{defw}
\rho_\lambda:=2^{\kappa |\lambda|},
\ee
satisfies 
\be
\sup_{x\in D} \sum_{\lambda} \rho_\lambda |\psi_\lambda(x)| \leq C M \sum_{l\geq 0} 2^{(\kappa-\alpha)l} <\infty,
\ee
which is \iref{firstcond} in Theorem \ref{maintheo}.

Note that if we order our wavelet basis from coarse to fine scale, we find for the 
resulting system $(\psi_j)_{j\geq 1}$ and sequence $(\rho_j)_{j\geq 1}$ the algebraic behaviour
\be
\|\psi_j\|_{L^\infty} \sim j^{-\alpha/d},
\ee
and
\be
\rho_j \sim j^{\kappa/d}.
\ee
We thus obtain from Theorem \ref{maintheo}
the following immediate corollary.

\begin{corollary}
\label{wavcor}
Let $(\psi_j)_{j\geq 1}$ be a wavelet basis with the normalization \iref{normalpha}.   If  
  $(\|\psi_j\|_{L^\infty})_{j\geq 1}$ belongs to $\ell^q(\N)$ 
then $(\|u_\nu\|_V)_{\nu\in \cF}$ belongs to $\ell^p(\cF)$  for all $p$ such that  $\frac{2p}{2-p}>q$.  In particular, the  best $n$-term Hermite approximations
converge in $L^2(U,V,\gamma)$ with rate $n^{-s}$ for all $s<\frac 1 q$.
\end{corollary}

\noindent
{\bf Proof:}  Since $(\|\psi_j\|_{L^\infty})_{j\geq 1}$ belongs to $\ell^q(\N)$, we know that  $q\alpha>d$.    Since for the given $p$, we have  $q^*:=\frac{2p}{2-p}>q$, we can     take  $ \kappa\in ]\frac{d}{q^*},\alpha[$ and $(\rho_\lambda)$ as defined in
\iref{defw}  for this $\kappa$.    Then, the conditions of  Theorem \ref{maintheo} are satisfied  for this sequence of weights and for  $q^*$.      
   Hence,  $(\|u_\nu\|_V)_{\nu\in \cF}$ belongs to $\ell^p(\cF)$.  This also gives that the approximation rate is $n^{-s}$ with $s=1/p-1/2=1/q^*$.  By adjusting $p$, we can take $q^*<q$ as close to $q$ as we wish, thereby establishing the rate of convergence 
   $n^{-s}$ for all $s<1/q$. \hfill $\Box$
\nl

Note that if we use sufficiently smooth wavelets, the decay property \iref{normalpha} is equivalent
to the property that the correlation function $C_b$ belongs to the Besov space $B^{\alpha}_{\infty}(L^\infty(D))$,
which coincides with the H\"older space $C^\alpha$ when $\alpha$ is non-integer,
and therefore $b$ is almost surely in the H\"older space $C^\beta$ for $\beta<\alpha/2$. Thus,
we also infer from Theorem \ref{maintheo} that if $C_b$ belongs to the
Besov space $B^{\alpha}_{\infty}(L^\infty(D))$ for some $\alpha>0$, best 
$n$-term Hermite approximations
converge in $L^2(U,V,\gamma)$ with rate $n^{-s}$ for all $s<\alpha/d$.

\subsection{Arbitrary supports}

We finally consider functions $\psi_j$ with
arbitrary supports, including the case of globally supported functions
such as the Fourier basis.  Let us assume that $(\|\psi_j\|_{L^\infty})_{j\geq 1}\in \ell^q(\N)$
for some $0<q<1$. We then find that the choice
\be
\rho_j:=\|\psi_j\|_{L^\infty}^{q-1},
\ee
obviously satisfies \iref{firstcond}. We also find that $(\rho_j^{-1})_{j\geq 1}$ belongs to $\ell^r(\N)$ for $r:=\frac{q}{1-q}$.   Therefore, applying Theorem \ref{maintheo}, we find that $(\|u_\nu\|_V)_{\nu\in\cF} \in \ell^p(\cF)$ when  $p$  satisfies
$r=\frac {2p}{2-p}$ or equivalently $\frac 1 q= \frac 1 p+\frac 1 2$.
Therefore, we obtain the following immediate corollary.

\begin{corollary}
\label{globcor}
Let $(\psi_j)_{j\geq 1}$ be a family of functions with arbitrary support, and let $0<p<2$. 
If $(\|\psi_j\|_{L^\infty})_{j\geq 1}$ belongs to $\ell^q(\N)$ with $\frac 1 q= \frac 1 p+\frac 1 2$,
then $(\|u_\nu\|_V)_{\nu\in \cF}$ belongs to $\ell^p(\cF)$. In particular, best $n$-term Hermite approximations
converge in $L^2(U,V,\gamma)$ with rate $n^{-s}$ for $s=\frac 1 q-1$.
\end{corollary}

Note that if $(j^\beta \|\psi_j\|_{L^\infty})_{j\geq 1}\in \ell^p(\N)$ for some $\beta>\frac 1 2$,
an application of H\"older's  inequality shows that $(\|\psi_j\|_{L^\infty})_{j\geq 1}$ belongs to 
$\ell^q(\N)$ with $\frac 1 q= \frac 1 p+\frac 1 2$. Therefore, the above corollary
represents an improvement
over the condition $(j \|\psi_j\|_{L^\infty})_{j\geq 1}\in \ell^p(\N)$ from \cite{HS}. 

\section{Non-optimality of the Karhunen-Lo\`eve representation}

The previous examples illustrate the role of the support properties
of the functions $(\psi_j)_{j\geq 1}$ when analyzing the convergence rate of 
best $n$-term Hermite approximations for the map $y\mapsto u(y)$.
In particular, they reveal that faster convergence rates can be obtained
in the case of locally supported functions.

We now discuss a concrete example which illustrates this phenomenon
for the approximation of a PDE with given lognormal coefficients, when 
we use two different representations of these coefficients. Here
we take
\be
D=]0,1[,
\ee
and therefore consider the equation
\be
-(au')'=f,\quad u(0)=u(1)=0.
\ee
We take $a=\exp(b)$ where $b$ is a Brownian bridge, that is, a Gaussian process
with covariance given by
\be
C_b(x,x')=\min\{x,x'\}-xx'.
\ee
There exists two simple explicit representations for this process.

\begin{enumerate}
\item
The Karhunen-Lo\`eve representation is determined by the eigenfunctions and eigenvalues of the covariance 
operator, which have the form
\be
\vp_j(x):=\sqrt 2\sin(\pi j x)\quad {\rm and} \quad \lambda_j:=\frac 1 {\pi^2 j^2},\quad j\geq 1,
\ee
so that after normalization, we obtain the representation 
\be
b(y)=\sum_{j\geq 1} y_j \psi_j,\quad  \psi_j(x):=\frac {\sqrt 2}{\pi j}\sin(\pi j x),
\ee
where the $y_j$ are i.i.d.\ standard Gaussian random variables.
\item
The Levy-Ciesielki representation uses the Schauder basis functions, which are the primitives of the
Haar functions, namely
\be
\psi_{l,k}(x)=2^{-l/2}\psi(2^l x-k), \quad  k=0,\dots,2^{l}-1,\; l\geq 0,
\ee
where $\psi(x):=\max\{0,1/2-|x-1/2|\}$. Then setting $\psi_j=\psi_{l,k}$ when $j=2^l+k$, we have again
\be
b(y)=\sum_{j\geq 1} y_j \psi_j,
\ee
where the $y_j$ are i.i.d.\ standard Gaussian random variables.
\end{enumerate}

If we use the Karhunen-Lo\`eve representation, the analysis carried in this paper
does not allow us to establish any $\ell^p$ summablity of $(\|u_\nu\|_V)_{\nu\in\cF}$ for $p<2$. Indeed,
due to the global nature of the functions $\psi_j$ we need to rely on Corollary
\ref{globcor} which requires that $(\|\psi_j\|_{L^\infty})_{j\geq 1}\in \ell^q(\N)$
for some $q<1$, which does not hold since this sequence is not $\ell^1$ summable.
Note however that $\ell^2$ summability is ensured since $\E(\|u(y)\|_V^2)<\infty$.

In contrast, when using the Levy-Ciesielski representation, which is essentially
of wavelet type with $\alpha=\frac 1 2$ in \iref{normalpha}, we can rely
on Corollary \ref{wavcor} which shows that $(\|u_\nu\|_V)_{\nu\in\cF}$ belongs
to $\ell^p(\cF)$ for all $1<p<2$. In particular, best $n$-term Hermite approximations
converge in $L^2(U,V,\gamma)$ with rate $n^{-s}$ for all $0<s<1$.

This example reveals that for a given lognormal process,
the Karhunen-Lo\`eve representation of the Gaussian process might
not be optimal in terms of the resulting convergence rates of the best $n$-term
Hermite approximation. One heuristic explanation of this fact is that the Karhunen-Lo\`eve
representation is optimal in a very specific sense: it minimizes the 
mean-square $L^2(D)$-error when truncating $b$ by the $J$ first terms
in its expansion. However, in the present setting of the elliptic diffusion equation, 
the relevant norm for approximating the functions $a$ and $b$ is 
not the $L^2$ norm, but rather the $L^\infty$ norm for which the Karhunen-Lo\`eve
representation has no particular optimality property.

\noindent
 Markus Bachmayr \\
 Sorbonne Universit\'es, UPMC Univ Paris 06, CNRS, UMR 7598, Laboratoire Jacques-Louis Lions, 4, place Jussieu 75005, Paris, France\\
 bachmayr@ljll.math.upmc.fr
 \vskip .1in
 \noindent
 Albert Cohen\\ 
Sorbonne Universit\'es, UPMC Univ Paris 06, CNRS, UMR 7598, Laboratoire Jacques-Louis Lions, 4, place Jussieu 75005, Paris, France \\
cohen@ljll.math.upmc.fr
\vskip .1in
 \noindent
Ronald DeVore\\
Department of Mathematics, Texas A\&M University,
College Station, TX 77840, USA\\
  rdevore@math.tamu.edu

 \vskip .1in
\noindent
 Giovanni  Migliorati \\
Sorbonne Universit\'es, UPMC Univ Paris 06, CNRS, UMR 7598, Laboratoire Jacques-Louis Lions, 4, place Jussieu 75005, Paris, France \\
migliorati@ljll.math.upmc.fr   

\end{document}